\documentclass[11pt,twoside]{article}
\usepackage{amsmath, amsthm, amscd, amsfonts, amssymb, graphicx, color}

\date{}

\begin{document}

\centerline{}

\centerline {\Large{\bf Generalized Riesz Representation Theorem in $n$-Hilbert space}}

\newcommand{\mvec}[1]{\mbox{\bfseries\itshape #1}}
\centerline{}
\centerline{\textbf{Prasenjit Ghosh}}
\centerline{Department of Pure Mathematics, University of Calcutta,}
\centerline{35, Ballygunge Circular Road, Kolkata, 700019, West Bengal, India}
\centerline{e-mail: prasenjitpuremath@gmail.com}
\centerline{}
\centerline{\textbf{T. K. Samanta}}
\centerline{Department of Mathematics, Uluberia College,}
\centerline{Uluberia, Howrah, 711315,  West Bengal, India}
\centerline{e-mail: mumpu$_{-}$tapas5@yahoo.co.in}

\newtheorem{Theorem}{\quad Theorem}[section]

\newtheorem{definition}[Theorem]{\quad Definition}

\newtheorem{theorem}[Theorem]{\quad Theorem}

\newtheorem{remark}[Theorem]{\quad Remark}

\newtheorem{corollary}[Theorem]{\quad Corollary}

\newtheorem{note}[Theorem]{\quad Note}

\newtheorem{lemma}[Theorem]{\quad Lemma}

\newtheorem{example}[Theorem]{\quad Example}

\newtheorem{result}[Theorem]{\quad Result}
\newtheorem{conclusion}[Theorem]{\quad Conclusion}

\newtheorem{proposition}[Theorem]{\quad Proposition}

\begin{abstract}
\textbf{\emph{In respect of \,$b$-linear functional, Riesz representation theorem in \,$n$-Hilbert space have been proved.\,We define \,$b$-sesquilinear functional in \,$n$-Hilbert space and establish the polarization identities.\,A generalized form of the Schwarz inequality in \,$n$-Hilbert space is being discussed.\,Finally, a generalized version of Riesz representation theorem with respect to \,$b$-sesquilinear functional in \,$n$-Hilbert space have been developed.}}
\end{abstract}
{\bf Keywords:}  \emph{Riesz representation theorem, sesquilinear functional, polarization  \\ \smallskip\hspace{2cm}identity, linear n-normed space, n-inner product space.}\\

\\

\section{Introduction}
 
\smallskip\hspace{.6 cm}The general form of a bounded linear functionals on various Banach spaces is quite difficult.\,However, in the special setting of Hilbert spaces, we get a representation theorem in terms of a fixed vector and the inner product, for any bounded linear functionals on the space.\,This theorem is known as the Riesz representation Theorem.\,From Riesz representation theorem it follows that the dual space \,$H^{\,\ast}$\, of a Hilbert space \,$H$\, is in one-to-one correspondence with the space \,$H$.\,This theorem is quite important in the theory of operators on Hilbert spaces.\,In particular, it refers  to the notion of Hilbert-adjoint operator of a bounded linear operator.\,This theorem is also used to present a general representation of sesquilinear functional on Hilbert space.

S.\,Gahler \cite{Gahler} introduced the notion of linear\;$2$-normed space.\;A geometric survey of the theory of linear\;$2$-normed space can be found in \cite{Freese}.\;The concept of $2$-Banach space is briefly discussed in \cite{White}.\;H.\,Gunawan and Mashadi \cite{Mashadi} developed the generalization of a linear\;$2$-normed space for \,$n \,\geq\, 2$.\,The concept of \,$2$-inner product space was first introduced by Diminnie et al. in 1970's \cite{Diminnie}.\;In 1989, A.\,Misiak \cite{Misiak} developed the generalization of a \,$2$-inner product space for \,$n \,\geq\, 2$.\,P. Ghosh and T. K. Samanta studied the Uniform Boundedness Principle and Hahn-Banach Theorem in linear\;$n$-normed space \cite{Prasenjit}.\,They also studied the reflexivity of linear\;$n$-normed space with respect to \,$b$-linear functional \cite{G}.   

In this paper,\;Riesz representation theorem for bounded\;$b$-linear functionals in case of \,$n$-Hilbert space is discussed.\;We present the notion of \,$b$-sesquilinear functional in \,$n$-Hilbert space and give some of its properties.\,The polarization identities associated with the \,$b$-sesquilinear functional in \,$n$-Hilbert space are given and a generalized form of the Schwarz inequality in \,$n$-Hilbert space is obtained.\,Finally, we present a general representation of bounded\;$b$-sesquilinear functional in \,$n$-Hilbert spaces.

\section{Preliminaries}

\begin{definition}\cite{Mashadi}
Let \,$H$\, be a linear space over the field \,$ \mathbb{K}$, where \,$ \mathbb{K} $\, is the real or complex numbers field with \,$\text{dim}\,H \,\geq\, n$, where \,$n$\, is a positive integer.\;A non-negative real valued function \,$\left \|\,\cdot \,,\, \cdots \,,\, \cdot \,\right \| \,:\, H^{\,n} \,\to\, \mathbb{R}$\, is called an n-norm on \,$X$\, provided for each \,$x,\, y,\, x_{\,1},\, x_{\,2},\, \cdots,\, x_{\,n} \,\in\, H$, 
\begin{description}
\item[$(i)$]$\left\|\,x_{\,1},\, x_{\,2},\, \cdots,\, x_{\,n}\,\right\| \,=\,0$\, if and only if \,$x_{\,1},\, \cdots,\, x_{\,n}$\, are linearly dependent,
\item[$(ii)$]$\left\|\,x_{\,1},\, x_{\,2},\, \cdots,\, x_{\,n}\,\right\|$\; is invariant under permutations of \,$x_{\,1},\, x_{\,2},\, \cdots,\, x_{\,n}$,
\item[$(iii)$]$\left\|\,\alpha\,x_{\,1},\, x_{\,2},\, \cdots,\, x_{\,n}\,\right\| \,=\, |\,\alpha\,|\, \left\|\,x_{\,1},\, x_{\,2},\, \cdots,\, x_{\,n}\,\right\|\; \;\;\forall \;\; \alpha \,\in\, \mathbb{K}$,
\item[$(iv)$]$\left\|\,x \,+\, y,\, x_{\,2},\, \cdots,\, x_{\,n}\,\right\| \,\leq\, \left\|\,x,\, x_{\,2},\, \cdots,\, x_{\,n}\,\right\| \,+\, \left\|\,y,\, x_{\,2},\, \cdots,\, x_{\,n}\,\right\|$.
\end{description}
A linear space \,$H$\, together with a \,$n$-norm \,$\left\|\,\cdot,\, \cdots,\, \cdot \,\right \|$\, on \,$H$\, is called a linear n-normed space.\;For particular value \,$n \,=\, 2$, the space \,$H$\, is said to be a linear 2-normed space \cite{Gahler}. 
\end{definition}

\begin{definition}\cite{Misiak}
Let \,$n \,\in\, \mathbb{N}$\; and \,$H$\, be a linear space of dimension greater than or equal to \,$n$\; over the field \,$\mathbb{K}$.\;A function \,$\left<\,\cdot,\, \cdot \,|\, \cdot,\, \cdots,\, \cdot\,\right> \,:\, H^{\, n\,+\, 1} \,\to\, \mathbb{K}$\, is called an \,$n$-inner product on \,$H$\, provided for all \,$x,\, y,\, x_{\,1},\, x_{\,2},\, \cdots,\, x_{\,n} \,\in\, H$,
\begin{description}
\item[$(i)$]\;\; $\left<\,x_{\,1},\, x_{\,1} \,|\, x_{\,2},\, \cdots,\, x_{\,n} \,\right> \,\geq\,  0$\, and \,$\left<\,x_{\,1},\, x_{\,1} \,|\, x_{\,2},\, \cdots,\, x_{\,n} \,\right> \,=\,  0$\, if and only if \,$x_{\,1},\, x_{\,2},\, \cdots,\, x_{\,n}$\, are linearly dependent,
\item[$(ii)$]\;\; $\left<\,x,\, y \,|\, x_{\,2},\, \cdots,\, x_{\,n} \,\right> \,=\, \left<\,x,\, y \,|\, x_{\,i_{\,2}},\, \cdots,\, x_{\,i_{\,n}}\,\right> $\, for every permutations \,$\left(\, i_{\,2},\, \cdots,\, i_{\,n} \,\right)$\, of \,$\left(\, 2,\, \cdots,\, n \,\right)$,
\item[$(iii)$]\;\; $\left<\,x,\, y \,|\, x_{\,2},\, \cdots,\, x_{\,n} \,\right> \,=\, \overline{\left<\,y,\, x \,|\, x_{\,2},\, \cdots,\, x_{\,n} \,\right> }$,
\item[$(iv)$]\;\; $\left<\,\alpha\,x,\, y \,|\, x_{\,2},\, \cdots,\, x_{\,n} \,\right> \,=\, \alpha \,\left<\,x,\, y \,|\, x_{\,2},\, \cdots,\, x_{\,n} \,\right> $, for \,$\alpha \,\in\, \mathbb{K}$,
\item[$(v)$]\;\; $\left<\,x \,+\, y,\, z \,|\, x_{\,2},\, \cdots,\, x_{\,n} \,\right> \,=\, \left<\,x,\, z \,|\, x_{\,2},\, \cdots,\, x_{\,n}\,\right> \,+\,  \left<\,y,\, z \,|\, x_{\,2},\, \cdots,\, x_{\,n} \,\right>$.
\end{description}
A linear space \,$H$\, together with n-inner product \,$\left<\,\cdot,\, \cdot \,|\, \cdot,\, \cdots,\, \cdot\,\right>$\, on \,$H$\, is called an n-inner product space.
\end{definition}

\begin{theorem}(Schwarz inequality)\cite{Misiak}\label{thn1}
Let \,$H$\, be a \,$n$-inner product space.\,Then 
\[\left|\,\left<\,x,\, y \,|\, x_{\,2},\, \cdots,\, x_{\,n}\,\right>\,\right| \,\leq\, \left\|\,x,\, x_{\,2},\, \cdots,\, x_{\,n}\,\right\|\, \left\|\,y,\, x_{\,2},\, \cdots,\, x_{\,n}\,\right\|\]
hold for all \,$x,\, y,\, x_{\,2},\, \cdots,\, x_{\,n} \,\in\, H$.
\end{theorem}

\begin{theorem}\cite{Misiak}
Let \,$H$\, be a \,$n$-inner product space.\,Then
\[\left \|\,x_{\,1},\, x_{\,2},\, \cdots,\, x_{\,n}\,\right\| \,=\, \sqrt{\left <\,x_{\,1},\, x_{\,1} \,|\, x_{\,2},\, \cdots,\, x_{\,n}\,\right>}\] defines a n-norm for which
\begin{align*}
&\|\,x \,+\, y,\, x_{\,2},\, \cdots,\, x_{\,n}\,\|^{\,2} \,+\, \|\,x \,-\, y,\, x_{\,2},\, \cdots,\, x_{\,n}\,\|^{\,2}\\
& \,=\, 2\, \left(\,\|\,x,\, x_{\,2},\, \cdots,\, x_{\,n}\,\|^{\,2} \,+\, \|\,y,\, x_{\,2},\, \cdots,\, x_{\,n}\,\|^{\,2} \,\right)
\end{align*} 
hold for all \,$x,\, y,\, x_{\,1},\, x_{\,2},\, \cdots,\, x_{\,n} \,\in\, H$.
\end{theorem}

\begin{definition}\cite{Mashadi}
Let \,$\left(\,H,\, \left\|\,\cdot,\, \cdots,\, \cdot \,\right\|\,\right)$\; be a linear n-normed space.\;A sequence \,$\{\,x_{\,k}\,\}$\; in \,$H$\, is said to convergent if there exists an \,$x \,\in\, H$\, such that 
\[\lim\limits_{k \to \infty}\,\left\|\,x_{\,k} \,-\, x,\, e_{\,2},\, \cdots,\, e_{\,n} \,\right\| \,=\, 0\]
for every \,$ e_{\,2},\, \cdots,\, e_{\,n} \,\in\, H$\, and it is called a Cauchy sequence if 
\[\lim\limits_{l,\, k \,\to\, \infty}\,\left \|\,x_{l} \,-\, x_{\,k},\, e_{\,2},\, \cdots,\, e_{\,n}\,\right\| \,=\, 0\]
for every \,$ e_{\,2},\, \cdots,\, e_{\,n} \,\in\, H$.\;The space \,$H$\, is said to be complete if every Cauchy sequence in this space is convergent in \,$H$.\;A n-inner product space is called n-Hilbert space if it is complete with respect to its induce norm.
\end{definition}

\begin{definition}\cite{Soenjaya}
We define the following open and closed ball in \,$H$: 
\[B_{\,\{\,e_{\,2},\, \cdots,\, e_{\,n}\,\}}\,(\,a,\, \delta\,) \,=\, \left\{\,x \,\in\, H \,:\, \left\|\,x \,-\, a,\, e_{\,2},\, \cdots,\, e_{\,n}\,\right\| \,<\, \delta \,\right\}\;\text{and}\]
\[B_{\,\{\,e_{\,2},\, \cdots,\, e_{\,n}\,\}}\,[\,a,\, \delta\,] \,=\, \left\{\,x \,\in\, H \,:\, \left\|\,x \,-\, a,\, e_{\,2},\, \cdots,\, e_{\,n}\,\right\| \,\leq\, \delta\,\right\},\hspace{.5cm}\]
where \,$a,\, e_{\,2},\, \cdots,\, e_{\,n} \,\in\, H$\, and \,$\delta$\, be a positive number.
\end{definition}

\begin{definition}\cite{Soenjaya}
A subset \,$G$\, of \,$H$\, is said to be open in \,$H$\, if for all \,$a \,\in\, G $, there exist \,$e_{\,2},\, \cdots,\, e_{\,n} \,\in\, H $\, and \, $\delta \,>\, 0 $\; such that \,$B_{\,\{\,e_{\,2},\, \cdots,\, e_{\,n}\,\}}\,(\,a,\, \delta\,) \,\subseteq\, G$.
\end{definition}

\begin{definition}\cite{Soenjaya}
Let \,$ A \,\subseteq\, H$.\;Then the closure of \,$A$\, is defined as 
\[\overline{A} \,=\, \left\{\, x \,\in\, H \;|\; \,\exists\, \;\{\,x_{\,k}\,\} \,\in\, A \;\;\textit{with}\;  \lim\limits_{k \,\to\, \infty} x_{\,k} \,=\, x \,\right\}.\]
The set \,$ A $\, is said to be closed if $ A \,=\, \overline{A}$. 
\end{definition}

\begin{definition}\cite{Prasenjit}\label{defn1}
Let \,$W$\, be a subspace of \,$H$\, and \,$b_{\,2},\, b_{\,3},\, \cdots,\, b_{\,n}$\, be fixed elements in \,$H$\, and \,$\left<\,b_{\,i}\,\right>$\, denote the subspaces of \,$H$\, generated by \,$b_{\,i}$, for \,$i \,=\, 2,\, 3,\, \cdots,\,n $.\;Then a map \,$T \,:\, W \,\times\,\left<\,b_{\,2}\,\right> \,\times\, \cdots \,\times\, \left<\,b_{\,n}\,\right> \,\to\, \mathbb{K}$\; is called a b-linear functional on \,$W \,\times\, \left<\,b_{\,2}\,\right> \,\times\, \cdots \,\times\, \left<\,b_{\,n}\,\right>$, if for every \,$x,\, y \,\in\, W$\, and \,$k \,\in\, \mathbb{K}$, the following conditions hold:
\begin{description}
\item[$(i)$] $T\,(\,x \,+\, y,\, b_{\,2},\, \cdots,\, b_{\,n}\,) \,=\, T\,(\,x,\, b_{\,2},\, \cdots,\, b_{\,n}\,) \,+\, T\,(\,y,\, b_{\,2},\, \cdots,\, b_{\,n}\,)$
\item[$(ii)$] $T\,(\,k\,x,\, b_{\,2},\, \cdots,\, b_{\,n}\,) \,=\, k\; T\,(\,x,\, b_{\,2},\, \cdots,\, b_{\,n}\,)$. 
\end{description}
A b-linear functional is said to be bounded if \,$\exists$\, a real number \,$M \,>\, 0$\, such that
\[\left|\,T\,(\,x,\, b_{\,2},\, \cdots,\, b_{\,n}\,)\,\right| \,\leq\, M\; \left\|\,x,\, b_{\,2},\, \cdots,\, b_{\,n}\,\right\|\; \;\forall\; x \,\in\, W.\]
The norm of the bounded b-linear functional \,$T$\, is defined by
\[\|\,T\,\| \,=\, \inf\,\left\{\,M \,>\, 0 \,:\, \left|\,T\,(\,x,\, b_{\,2},\, \cdots,\, b_{\,n}\,)\,\right| \,\leq\, M\; \left\|\,x,\, b_{\,2},\, \cdots,\, b_{\,n}\,\right\| \;\forall\; x \,\in\, W\,\right\}.\]
The norm of \,$T$\, can be expressed by any one of the following equivalent formula:
\begin{description}
\item[$(i)$]$\|\,T\,\| \,=\, \sup\,\left\{\,\left|\,T\,(\,x,\, b_{\,2},\, \cdots,\, b_{\,n}\,)\,\right| \;:\; \left\|\,x,\, b_{\,2},\, \cdots,\, b_{\,n}\,\right\| \,\leq\, 1\,\right\}$.
\item[$(ii)$]$\|\,T\,\| \,=\, \sup\,\left\{\,\left|\,T\,(\,x,\, b_{\,2},\, \cdots,\, b_{\,n}\,)\,\right| \;:\; \left\|\,x,\, b_{\,2},\, \cdots,\, b_{\,n}\,\right\| \,=\, 1\,\right\}$.
\item[$(iii)$]$ \|\,T\,\| \,=\, \sup\,\left \{\,\dfrac{\left|\,T\,(\,x,\, b_{\,2},\, \cdots,\, b_{\,n}\,)\,\right|}{\left\|\,x,\, b_{\,2},\, \cdots,\, b_{\,n}\,\right\|} \;:\; \left\|\,x,\, b_{\,2},\, \cdots,\, b_{\,n}\,\right\| \,\neq\, 0\,\right \}$. 
\end{description}
Also, we have \,$\left|\,T\,(\,x,\, b_{\,2},\, \cdots,\, b_{\,n}\,)\,\right| \,\leq\, \|\,T\,\|\, \left\|\,x,\, b_{\,2},\, \cdots,\, b_{\,n}\,\right\|\, \;\forall\; x \,\in\, W$.
\end{definition}
Let \,$H_{F}^{\,\ast}$\, denotes the Banach space of all bounded\;$b$-linear functional defined on \,$H \,\times\, \left<\,b_{\,2}\,\right> \,\times \cdots \,\times\, \left<\,b_{\,n}\,\right>$\, with respect to the above norm.

\section{Riesz Representation Theorem in $n$-Hilbert space}

\smallskip\hspace{.6 cm}In this section, we explore a relationship between the vectors in the \,$n$-Hilbert space \,$H$\, and the bounded\;$b$-linear functionals defined on  \,$H \,\times\,\left<\,b_{\,2}\,\right> \,\times\, \cdots \,\times\, \left<\,b_{\,n}\,\right>$. 

\begin{definition}
Let \,$S$\, be a subset of a n-Hilbert space \,$H$.\,Two elements \,$x$\, and \,$y$\, of \,$H$\, are said to be b-orthogonal if \,$\left<\,x,\,y \,|\, b_{\,2},\, \cdots,\, b_{\,n}\,\right> \,=\, 0$.\,In symbol, we write \,$x \,\perp\, y$.\,If \,$x$\, is b-orthogonal to every element of \,$S$, then we say that \,$x$\, is b-orthogonal to \,$S$\, and in symbol we write \,$x \,\perp\, S$.
\end{definition} 

\begin{definition}
Let \,$S \,\subseteq\, H$.\,Then the set of all elements of \,$H$, b-orthogonal to \,$S$\, is called the b-orthogonal complement of \,$S$\, and is denoted by \,$S^{\,\perp}$.
\end{definition}

\begin{theorem}\label{th2}
Let \,$H$\, be a \,$n$-Hilbert space.\,Then \,$T$\, is a bounded\;$b$-linear functional defined on \,$H \,\times\, \left<\,b_{\,2}\,\right> \,\times\, \cdots \,\times\, \left<\,b_{\,n}\,\right>$\, if and only if there exists a unique element \,$z$\, in \,$H$\, with \,$\left\{\,z,\, b_{\,2},\, \cdots,\, b_{\,n}\,\right\}$\, is linearly independent such that
\begin{equation}\label{eq2} 
T\,(\,x,\, b_{\,2},\, \cdots,\, b_{\,n}\,) \,=\, \left<\,x,\, z \,|\, b_{\,2},\, \cdots,\, b_{\,n}\,\right>, \;\text{for all \,$x \,\in\, H$}
\end{equation}
and moreover \,$\|\,T\,\| \,=\, \left\|\,z,\,b_{\,2},\, \cdots,\, b_{\,n}\,\right\|$.
\end{theorem} 

\begin{proof}
Let \,$z$\, be any fixed element in \,$H$\, and define a functional \,$T$\, by
\[T\,(\,x,\, b_{\,2},\, \cdots,\, b_{\,n}\,) \,=\, \left<\,x,\, z \,|\, b_{\,2},\, \cdots,\, b_{\,n}\,\right>, \;\text{for all \,$x \,\in\, H$}.\]Then
\begin{description}
\item[$(i)$]$T$\, is a \,$b$-linear functional defined on \,$H \,\times\, \left<\,b_{\,2}\,\right> \,\times\, \cdots \,\times\, \left<\,b_{\,n}\,\right>$;
\begin{align*}
T\,(\,x \,+\, y,\, b_{\,2},\, \cdots,\, b_{\,n}\,)&\,=\, \left<\,x \,+\, y,\, z \,|\, b_{\,2},\, \cdots,\, b_{\,n}\,\right>\\
&\,=\, \left<\,x,\, z \,|\, b_{\,2},\, \cdots,\, b_{\,n}\,\right> \,+\, \left<\,y,\, z \,|\, b_{\,2},\, \cdots,\, b_{\,n}\,\right>\\
&=\, T\,(\,x,\, b_{\,2},\, \cdots,\, b_{\,n}\,) \,+\, T\,(\,y,\, b_{\,2},\, \cdots,\, b_{\,n}\,),\;\text{and}\\
T\,(\,k\,x,\, b_{\,2},\, \cdots,\, b_{\,n}\,) &\,=\, \left<\,k\,x,\, z \,|\, b_{\,2},\, \cdots,\, b_{\,n}\,\right> \,=\, k\,\left<\,x,\, z \,|\, b_{\,2},\, \cdots,\, b_{\,n}\,\right>\hspace{1cm}\\
&\,=\, k\,T\,(\,x,\, b_{\,2},\, \cdots,\, b_{\,n}\,),\; \;\text{for all}\; \,x,\,y \,\in\, H\; \;\text{and}\; \,k \,\in\, \mathbb{K}.
\end{align*}
\item[$(ii)$]$T$\, is bounded;
\begin{align*}
\left|\,T\,(\,x,\, b_{\,2},\, \cdots,\, b_{\,n}\,)\,\right| &\,=\, \left|\,\left<\,x,\, z \,|\, b_{\,2},\, \cdots,\, b_{\,n}\,\right>\,\right|\\
&\leq\, \left\|\,x,\,b_{\,2},\, \cdots,\, b_{\,n}\,\right\|\,\left\|\,z,\,b_{\,2},\, \cdots,\, b_{\,n}\,\right\|,
\end{align*}
for all \,$x \,\in\, H$. 
\end{description}
Since \,$z$\, is fixed, the above calculation shows that \,$T$\, is bounded and moreover \,$\|\,T\,\| \,\leq\, \left\|\,z,\,b_{\,2},\, \cdots,\, b_{\,n}\,\right\|$.\,On the other hand, if \,$z \,\neq\, \theta$\, with \,$\left\{\,z,\, b_{\,2},\, \cdots,\, b_{\,n}\,\right\}$\, is linearly independent then
\begin{align*}
\|\,T\,\| & \,=\, \sup\,\left\{\,\left|\,T\,(\,x,\, b_{\,2},\, \cdots,\, b_{\,n}\,)\,\right| \;:\; x \,\in\, H, \left\|\,x,\, b_{\,2},\, \cdots,\, b_{\,n}\,\right\| \,\leq\, 1\,\right\}\\
&=\, \sup\,\left\{\,\left|\,\left<\,x,\, z \,|\, b_{\,2},\, \cdots,\, b_{\,n}\,\right>\,\right| \;:\; x \,\in\, H, \left\|\,x,\, b_{\,2},\, \cdots,\, b_{\,n}\,\right\| \,\leq\, 1\,\right\}\\
&\geq\, \left<\,\dfrac{z}{\left\|\,z,\, b_{\,2},\, \cdots,\, b_{\,n}\,\right\|},\, z \,|\, b_{\,2},\, \cdots,\, b_{\,n}\,\right> \,=\, \left\|\,z,\, b_{\,2},\, \cdots,\, b_{\,n}\,\right\|.
\end{align*}
In case \,$z \,=\, \theta$\, or \,$\left\{\,z,\, b_{\,2},\, \cdots,\, b_{\,n}\,\right\}$\, is linearly dependent \,$\|\,T\,\| \,\geq\, \left\|\,z,\,b_{\,2},\, \cdots,\, b_{\,n}\,\right\|$\, is obvious.\,Hence \,$\|\,T\,\| \,=\, \left\|\,z,\,b_{\,2},\, \cdots,\, b_{\,n}\,\right\|$.\\

Conversely, suppose that \,$T$\, is a bounded\;$b$-linear functional defined on \,$H \,\times\, \left<\,b_{\,2}\,\right> \,\times\, \cdots \,\times\, \left<\,b_{\,n}\,\right>$.\,If \,$T \,=\, 0$\, then the proof hold if we take \,$z \,=\, \theta$\, or \,$\left\{\,z,\, b_{\,2},\, \cdots,\, b_{\,n}\,\right\}$\, is linearly dependent.\,Let \,$T \,\neq\, 0$.\,Since \,$T$\, is bounded\;$b$-linear functional, the null space \,$\mathcal{N}\,(\,T\,)$\, is a closed subspace of \,$H$.\,Because \,$T \,\neq\, 0$, it follows that \,$\mathcal{N}\,(\,T\,) \,\neq\, H$.\,But, by the projection theorem \,$H \,=\, \mathcal{N} \,\oplus\, \mathcal{N}^{\,\perp}$\, with \,$\mathcal{N} \,=\, \mathcal{N}\,(\,T\,)$\, and therefore \,$\mathcal{N}^{\,\perp} \,\neq\, \{\,\theta\,\}$.\,This implies that there exists an element \,$z_{\,0} \,\in\, \mathcal{N}^{\,\perp}$\, such that \,$z_{\,0} \,\neq\, \theta$.\,Consider the set
\[S \,=\, \left\{\,v \,=\, z_{\,0}\,T\,(\,x,\, b_{\,2},\, \cdots,\, b_{\,n}\,) \,-\, x\,T\,(\,z_{\,0},\, b_{\,2},\, \cdots,\, b_{\,n}\,) \,:\, x \,\in\, H\,\right\}.\]
Then \,$S \,\subset\, \mathcal{N}$, since
\begin{align*}
&T\,(\,v,\, b_{\,2},\, \cdots,\, b_{\,n}\,)\\
 & \,=\, T\,\left[\,z_{\,0}\,T\,(\,x,\, b_{\,2},\, \cdots,\, b_{\,n}\,) \,-\, x\,T\,(\,z_{\,0},\, b_{\,2},\, \cdots,\, b_{\,n}\,) \,,\, b_{\,2},\, \cdots,\, b_{\,n}\,\right]\\
&=\, T\,(\,z_{\,0},\, b_{\,2},\, \cdots,\, b_{\,n}\,)\,T\,(\,x,\, b_{\,2},\, \cdots,\, b_{\,n}\,) \,-\, \,T\,(\,x,\, b_{\,2},\, \cdots,\, b_{\,n}\,)\,T\,(\,z_{\,0},\, b_{\,2},\, \cdots,\, b_{\,n}\,)\\
&=\, 0, \;\text{for all}\; x \,\in\, H.
\end{align*} 
Therefore, \,$z_{\,0} \,\perp\, S$.\,This gives
\begin{align*}
&\left<\,z_{\,0}\,T\,(\,x,\, b_{\,2},\, \cdots,\, b_{\,n}\,) \,-\, x\,T\,(\,z_{\,0},\, b_{\,2},\, \cdots,\, b_{\,n}\,),\, z_{\,0} \,|\, b_{\,2},\, \cdots,\, b_{\,n}\,\right> \,=\, 0\\
&\Rightarrow\,T\,(\,x,\, b_{\,2},\, \cdots,\, b_{\,n}\,)\,\left\|\,z_{\,0},\,b_{\,2},\, \cdots,\, b_{\,n}\,\right\|^{\,2} \,=\, T\,(\,z_{\,0},\, b_{\,2},\, \cdots,\, b_{\,n}\,)\,\left<\,x,\, z_{\,0} \,|\, b_{\,2},\, \cdots,\, b_{\,n}\,\right>.
\end{align*}
This implies that
\begin{align*}
T\,(\,x,\, b_{\,2},\, \cdots,\, b_{\,n}\,) &\,=\, \dfrac{T\,(\,z_{\,0},\, b_{\,2},\, \cdots,\, b_{\,n}\,)}{\left\|\,z_{\,0},\,b_{\,2},\, \cdots,\, b_{\,n}\,\right\|^{\,2}}\,\left<\,x,\, z_{\,0} \,|\, b_{\,2},\, \cdots,\, b_{\,n}\,\right>\\
& \,=\,\left<\,x,\, z \,|\, b_{\,2},\, \cdots,\, b_{\,n}\,\right>\; \;\text{for all}\; \,x \,\in\, H, 
\end{align*}
where \[z \,=\, \dfrac{\overline{T\,(\,z_{\,0},\, b_{\,2},\, \cdots,\, b_{\,n}\,)}\,z_{\,0}}{\left\|\,z_{\,0},\,b_{\,2},\, \cdots,\, b_{\,n}\,\right\|^{\,2}} \,\in\, H.\]  
This proves the existence of \,$z$.\;Let \,$z_{\,1},\, z_{\,2} \,\in\, H$\, with \,$\left\{\,z_{\,1} \,-\, z_{\,2},\,b_{\,2},\, \cdots,\, b_{\,n}\,\right\}$\, be linearly independent such that
\begin{align*}
T\,(\,x,\, b_{\,2},\, \cdots,\, b_{\,n}\,)& \,=\, \left<\,x,\, z_{\,1} \,|\, b_{\,2},\, \cdots,\, b_{\,n}\,\right>\\
& \,=\, \left<\,x,\, z_{\,2} \,|\, b_{\,2},\, \cdots,\, b_{\,n}\,\right>, \;\text{for all}\; x \,\in\, H.
\end{align*}  
Then 
\[\left<\,x,\, z_{\,1} \,-\, z_{\,2} \,|\, b_{\,2},\, \cdots,\, b_{\,n}\,\right> \,=\, 0\; \;\text{for all}\; \,x \,\in\, H.\]\,In particular, for \,$x \,=\, z_{\,1} \,-\, z_{\,2}$,
\begin{align*}
&\left<\,z_{\,1} \,-\, z_{\,2},\, z_{\,1} \,-\, z_{\,2} \,|\, b_{\,2},\, \cdots,\, b_{\,n}\,\right> \,=\, 0\\
&\;\Rightarrow\, \left\|\,z_{\,1} \,-\, z_{\,2},\, b_{\,2},\, \cdots,\, b_{\,n}\,\right\|^{\,2} \,=\, 0.
\end{align*}
This implies that \,$z_{\,1} \,-\, z_{\,2} \,=\, 0\;\Rightarrow\, z_{\,1} \,=\, z_{\,2}$.\,This proves the uniqueness of \,$z$.\\
Now, for \,$x \,=\, z$\, in (\ref{eq2}), we have 
\[T\,(\,z,\, b_{\,2},\, \cdots,\, b_{\,n}\,) \,=\, \left<\,z,\, z \,|\, b_{\,2},\, \cdots,\, b_{\,n}\,\right> \,=\, \left\|\,z,\, b_{\,2},\, \cdots,\, b_{\,n}\,\right\|^{\,2}.\]
But, \,$T$\, is bounded, we get
\[\left\|\,z,\, b_{\,2},\, \cdots,\, b_{\,n}\,\right\|^{\,2} \,=\, \left|\,T\,(\,z,\, b_{\,2},\, \cdots,\, b_{\,n}\,)\,\right| \,\leq\, \|\,T\,\|\,\left\|\,z,\, b_{\,2},\, \cdots,\, b_{\,n}\,\right\|,\] and so that \,$\left\|\,z,\, b_{\,2},\, \cdots,\, b_{\,n}\,\right\| \,\leq\, \|\,T\,\|$.\,On the other hand, using Schwartz inequality, 
\begin{align*}
\left|\,T\,(\,x,\, b_{\,2},\, \cdots,\, b_{\,n}\,)\,\right| &\,=\, \left|\,\left<\,x,\, z \,|\, b_{\,2},\, \cdots,\, b_{\,n}\,\right>\,\right|\\
& \,\leq\, \left\|\,x,\, b_{\,2},\, \cdots,\, b_{\,n}\,\right\|\,\left\|\,z,\, b_{\,2},\, \cdots,\, b_{\,n}\,\right\|,
\end{align*}
\[\Rightarrow\,\|\,T\,\| \,=\, \sup\limits_{\left\|\,x,\, b_{\,2},\, \cdots,\, b_{\,n}\,\right\| \,\leq\, 1}\,\left|\,\left<\,x,\, z \,|\, b_{\,2},\, \cdots,\, b_{\,n}\,\right>\,\right| \,\leq\, \left\|\,z,\, b_{\,2},\, \cdots,\, b_{\,n}\,\right\|.\hspace{2.5cm}\]
Hence, \,$\|\,T\,\| \,=\, \left\|\,z,\,b_{\,2},\, \cdots,\, b_{\,n}\,\right\|$.\,This completes the proof. 
\end{proof} 
 
\section{$b$-sesquilinear functional in $n$-Hilbert space}

\smallskip\hspace{.6 cm}In this section, we introduce the concept of bounded\;$b$-sesquilinear functional and discuss some of its properties.\,Finally, we present a general representation of \,$b$-sesquilinear  in \,$n$-Hilbert spaces. 

\begin{definition}
Let \,$H$\, be a linear space over the field \,$\mathbb{K}$.\,A \,$b$-sesquilinear functional \,$T$\, defined on \,$H \,\times\, H \,\times\,\left<\,b_{\,2}\,\right> \,\times\, \cdots \,\times\, \left<\,b_{\,n}\,\right>$\, is a mapping 
\[T \,:\, H \,\times\, H \,\times\,\left<\,b_{\,2}\,\right> \,\times\, \cdots \,\times\, \left<\,b_{\,n}\,\right> \,\to\, \mathbb{K}\] which satisfies the following conditions:
\begin{description}
\item[$(i)$]$T\,(\,x \,+\, y,\, z,\, b_{\,2},\, \cdots,\, b_{\,n}\,) \,=\, T\,(\,x,\, z,\, b_{\,2},\, \cdots,\, b_{\,n}\,) \,+\, T\,(\,y,\, z,\, b_{\,2},\, \cdots,\, b_{\,n}\,)$,
\item[$(ii)$] $T\,(\,\alpha\,x,\, y,\, b_{\,2},\, \cdots,\, b_{\,n}\,) \,=\, \alpha\,T\,(\,x,\, y,\, b_{\,2},\, \cdots,\, b_{\,n}\,)$, 
\item[$(iii)$]$T\,(\,x,\, y \,+\, z,\, b_{\,2},\, \cdots,\, b_{\,n}\,) \,=\, T\,(\,x,\, y,\, b_{\,2},\, \cdots,\, b_{\,n}\,) \,+\, T\,(\,x,\, z,\, b_{\,2},\, \cdots,\, b_{\,n}\,)$,
\item[$(iv)$] $T\,(\,x,\, \beta\,y,\, b_{\,2},\, \cdots,\, b_{\,n}\,) \,=\, \overline{\beta}\,T\,(\,x,\, y,\, b_{\,2},\, \cdots,\, b_{\,n}\,)$, 
\end{description}
for all \,$x,\, y,\, z \,\in\, H$\, and \,$\alpha,\, \beta \,\in\, \mathbb{K}$.
\end{definition} 
 
\begin{example}
If \,$H$\, is an \,$n$-inner product space and if we define a function  \,$T \,:\, H \,\times\, H \,\times\,\left<\,b_{\,2}\,\right> \,\times\, \cdots \,\times\, \left<\,b_{\,n}\,\right> \,\to\, \mathbb{K}$\, by
\[T\,(\,x,\, y,\, b_{\,2},\, \cdots,\, b_{\,n}\,) \,=\, \left<\,x,\, y \,|\, b_{\,2},\, \cdots,\, b_{\,n}\,\right>, \;\text{for all}\; \,x,\, y \,\in\, H.\]
Then \,$T$\, is a b-sesquilinear functional.
\end{example} 

\begin{example}
Let \,$T \,:\, H \,\times\, H \,\times\,\left<\,b_{\,2}\,\right> \,\times\, \cdots \,\times\, \left<\,b_{\,n}\,\right> \,\to\, \mathbb{K}$\, be a b-sesquilinear functional.\,We define a functional \,$U$\, as 
\[U\,(\,x,\, y,\, b_{\,2},\, \cdots,\, b_{\,n}\,) \,=\, \overline{\,T\,(\,y,\, x,\, b_{\,2},\, \cdots,\, b_{\,n}\,)}, \;\text{for all}\; \,x,\, y \,\in\, H.\]
Then \,$U$\, is a b-sesquilinear functional. 
\end{example}
 
\begin{example}
If \,$H$\, is a n-inner product space and \,$A \,:\, H \,\to\, H$\, is a linear operator then 
\[T\,(\,x,\, y,\, b_{\,2},\, \cdots,\, b_{\,n}\,) \,=\, \left<\,A\,x,\, y \,|\, b_{\,2},\, \cdots,\, b_{\,n}\,\right>,\; \;\text{for all}\; \,x,\, y \,\in\, H\] is a \,$b$-sesquilinear functional. 
\end{example}

\begin{definition}
Let \,$T \,:\, H \,\times\, H \,\times\,\left<\,b_{\,2}\,\right> \,\times\, \cdots \,\times\, \left<\,b_{\,n}\,\right> \,\to\, \mathbb{K}$\, be a \,$b$-sesquilinear functional.\,Then
\begin{description}
\item[$(i)$] If \,$T\,(\,x,\, y,\, b_{\,2},\, \cdots,\, b_{\,n}\,) \,=\, \overline{\,T\,(\,y,\, x,\, b_{\,2},\, \cdots,\, b_{\,n}\,)}$, for all \,$x,\, y \,\in\, H$, then \,$T$\, is called a symmetric \,$b$-sesquilinear functional.
\item[$(ii)$] If \,$T\,(\,x,\, x,\, b_{\,2},\, \cdots,\, b_{\,n}\,) \,\geq\, 0$, for all \,$x \,\in\, H$, then \,$T$\, is called positive.
\item[$(iii)$]The map \,$T^{\,\prime} \,:\, H \,\times\,\left<\,b_{\,2}\,\right> \,\times\, \cdots \,\times\, \left<\,b_{\,n}\,\right> \,\to\, \mathbb{R}$\, defined by
\[T^{\,\prime}\,(\,x,\, b_{\,2},\, \cdots,\, b_{\,n}\,) \,=\, T\,(\,x,\, x,\, b_{\,2},\, \cdots,\, b_{\,n}\,),\; \;x \,\in\, H\]  is called the quadratic form associated with the \,$b$-sesquilinear functional \,$T$.  
\end{description}
\end{definition} 

\begin{theorem}\label{lem2}(Polarization identities)
If \,$T$\, is a symmetric \,$b$-sesquilinear functional defined on \,$H \,\times\, H \,\times\,\left<\,b_{\,2}\,\right> \,\times\, \cdots \,\times\, \left<\,b_{\,n}\,\right>$\, and \,$T^{\,\prime}$\, is the associated quadratic form.\,Then we have the followings:
\begin{description}
\item[$(i)$]If \,$\mathbb{K} \,=\, \mathbb{R}$, then for all \,$x,\, y \,\in\, H$, we have
\[T\,(\,x,\, y,\, b_{\,2},\, \cdots,\, b_{\,n}\,) \,=\, \dfrac{1}{4}\,\left[\,T^{\,\prime}\,(\,x \,+\, y,\, b_{\,2},\, \cdots,\, b_{\,n}\,) \,-\, T^{\,\prime}\,(\,x \,-\, y,\, b_{\,2},\, \cdots,\, b_{\,n}\,)\,\right].\]
\item[$(ii)$]If \,$\mathbb{K} \,=\, \mathbb{C}$, then for all \,$x,\, y \,\in\, H$, we have
\begin{align*}
&T\,(\,x,\, y,\, b_{\,2},\, \cdots,\, b_{\,n}\,)\\
&=\, \dfrac{1}{4}\,\left[\,T^{\,\prime}\,(\,x \,+\, y,\, b_{\,2},\, \cdots,\, b_{\,n}\,) \,-\, T^{\,\prime}\,(\,x \,-\, y,\, b_{\,2},\, \cdots,\, b_{\,n}\,)\,\right] \,+\\
&+\,\dfrac{1}{4}\,\left[\,i\,T^{\,\prime}\,(\,x \,+\, i\,y,\, b_{\,2},\, \cdots,\, b_{\,n}\,) \,-\, i\,T^{\,\prime}\,(\,x \,-\, i\,y,\, b_{\,2},\, \cdots,\, b_{\,n}\,)\,\right].
\end{align*}
\end{description}  
\end{theorem} 

\begin{proof}
For every \,$x,\, y \,\in\, H$, we have
\begin{align*}
&T^{\,\prime}\,(\,x \,+\, y,\, b_{\,2},\, \cdots,\, b_{\,n}\,) \,-\, T^{\,\prime}\,(\,x \,-\, y,\, b_{\,2},\, \cdots,\, b_{\,n}\,)\\
&=\, T\,(\,x \,+\, y,\, x \,+\, y,\, b_{\,2},\, \cdots,\, b_{\,n}\,) \,-\, T\,(\,x \,-\, y,\, x \,-\, y,\, b_{\,2},\, \cdots,\, b_{\,n}\,)\\ 
&=\, 2\,\left[\,T\,(\,x,\, y,\, b_{\,2},\, \cdots,\, b_{\,n}\,) \,+\, \overline{T\,(\,x,\, y,\, b_{\,2},\, \cdots,\, b_{\,n}\,)}\,\right]\\
& \,=\, 4\,\text{Re}\,T\,(\,x,\, y,\, b_{\,2},\, \cdots,\, b_{\,n}\,).
\end{align*}
If the scalar field is the set of complex numbers, then for every \,$x,\, y \,\in\, H$, we have
\begin{align*}
&T^{\,\prime}\,(\,x \,+\, i\,y,\, b_{\,2},\, \cdots,\, b_{\,n}\,) \,-\, T^{\,\prime}\,(\,x \,-\, i\,y,\, b_{\,2},\, \cdots,\, b_{\,n}\,)\\
& \,=\, 4\,\text{Im}\,T\,(\,x,\, y,\, b_{\,2},\, \cdots,\, b_{\,n}\,).
\end{align*}
Thus, if \,$\mathbb{K} \,=\, \mathbb{C}$, then adding the above two equalities, we get \,$(ii)$\, and if \,$\mathbb{K} \,=\, \mathbb{R}$, then we have
\[T\,(\,x,\, y,\, b_{\,2},\, \cdots,\, b_{\,n}\,) \,=\, \dfrac{1}{4}\,\left(\,T^{\,\prime}\,(\,x \,+\, y,\, b_{\,2},\, \cdots,\, b_{\,n}\,) \,-\, T^{\,\prime}\,(\,x \,-\, y,\, b_{\,2},\, \cdots,\, b_{\,n}\,)\,\right).\]
The relations \,$(i)$\, and \,$(ii)$\, are called polarization identities associated with the \,$b$-sesquilinear functional. 
\end{proof}

\begin{theorem}\label{th2.1}
Let \,$T$\, be a \,$b$-sesquilinear functional and \,$T^{\,\prime}$\, be its associated quadratic form.\,Then \,$T$\, is symmetric if and only if \,$T^{\,\prime}$\, is real-valued. 
\end{theorem} 
 
\begin{proof}
Suppose \,$T$\, is symmetric and so 
\[T\,(\,x,\, y,\, b_{\,2},\, \cdots,\, b_{\,n}\,) \,=\, \overline{T\,(\,y,\, x,\, b_{\,2},\, \cdots,\, b_{\,n}\,)},\; \;\text{for all}\; \,x,\,y \,\in\, H.\]
Then, for all \,$x \,\in\, H$, we have
\begin{align*}
T^{\,\prime}\,(\,x,\, b_{\,2},\, \cdots,\, b_{\,n}\,) &\,=\, T\,(\,x,\, x,\, b_{\,2},\, \cdots,\, b_{\,n}\,) \,=\, \overline{T\,(\,x,\, x,\, b_{\,2},\, \cdots,\, b_{\,n}\,)}\\
& \,=\, \overline{T^{\,\prime}\,(\,x,\, b_{\,2},\, \cdots,\, b_{\,n}\,)}.
\end{align*}
This shows that \,$T^{\,\prime}$\, is real-valued.\\

Conversely, suppose that \,$T^{\,\prime}$\, is real-valued and let \[U\,(\,x,\, y,\, b_{\,2},\, \cdots,\, b_{\,n}\,) \,=\, \overline{T\,(\,y,\, x,\, b_{\,2},\, \cdots,\, b_{\,n}\,)},\; \;\text{for all}\; \,x,\,y \,\in\, H.\]
Then, for all \,$x \,\in\, H$, we have   
\begin{align*}
 U^{\,\prime}\,(\,x,\, b_{\,2},\, \cdots,\, b_{\,n}\,) &\,=\, U\,(\,x,\, x,\, b_{\,2},\, \cdots,\, b_{\,n}\,) \,=\, \overline{T\,(\,x,\, x,\, b_{\,2},\, \cdots,\, b_{\,n}\,)}\\
&=\,  \overline{T^{\,\prime}\,(\,x,\, b_{\,2},\, \cdots,\, b_{\,n}\,)} \,=\, T^{\,\prime}\,(\,x,\, b_{\,2},\, \cdots,\, b_{\,n}\,).
\end{align*}
Now, using Theorem (\ref{lem2}), it follows that \,$U \,=\, T$\, i.\,e.,
\[T\,(\,x,\, y,\, b_{\,2},\, \cdots,\, b_{\,n}\,) \,=\, \overline{T\,(\,y,\, x,\, b_{\,2},\, \cdots,\, b_{\,n}\,)},\; \;\text{for all}\; \,x,\,y \,\in\, H.\]So, \,$T$\, is symmetric.\,This proves the theorem.  
\end{proof} 

\begin{definition}
Let \,$T \,:\, H \,\times\, H \,\times\,\left<\,b_{\,2}\,\right> \,\times\, \cdots \,\times\, \left<\,b_{\,n}\,\right> \,\to\, \mathbb{K}$\, be a \,$b$-sesquilinear functional and \,$T^{\,\prime}$\, be its associated quadratic form.\,Then
\begin{description}
\item[$(i)$]\,$T$\, is said to be bounded if there exists \,$M \,>\, 0$\, such that
\[\left|\,T\,(\,x,\, y,\, b_{\,2},\, \cdots,\, b_{\,n}\,)\,\right| \,\leq\, M\,\left\|\,x,\,b_{\,2},\, \cdots,\, b_{\,n}\,\right\|\,\left\|\,y,\,b_{\,2},\, \cdots,\, b_{\,n}\,\right\|,\; \;\forall\; x,\, y \,\in\, H.\]
The infimum of all such \,$M$, is called the norm of \,$T$\, and is denoted by \,$\|\,T\,\|$.
\item[$(ii)$]\,$T^{\,\prime}$\, is said to be bounded if there exists \,$M \,>\, 0$\, such that
\[\left|\,T^{\,\prime}\,(\,x,\, b_{\,2},\, \cdots,\, b_{\,n}\,)\,\right| \,\leq\, M\; \left\|\,x,\, b_{\,2},\, \cdots,\, b_{\,n}\,\right\|^{\,2},\; \;\text{for all}\; \,x \,\in\, H.\]
The infimum of all such \,$M$, is called the norm of \,$T^{\,\prime}$\, and is denoted by \,$\|\,T^{\,\prime}\,\|$.
\end{description}
\end{definition}

\begin{remark}\label{note2}
According to the definition (\ref{defn1}), we can write
\begin{description}
\item[$(i)$]$\left|\,T\,(\,x,\, y,\, b_{\,2},\, \cdots,\, b_{\,n}\,)\,\right| \,\leq\, \|\,T\,\|\,\left\|\,x,\,b_{\,2},\, \cdots,\, b_{\,n}\,\right\|\,\left\|\,y,\,b_{\,2},\, \cdots,\, b_{\,n}\,\right\|,\\ \;\text{for all}\; \,x,\, y \,\in\, H$.
\item[$(ii)$]$\|\,T\,\| \,=\, \sup\limits_{\left\|\,x,\, b_{\,2},\, \cdots,\, b_{\,n}\,\right\| \,=\, 1 \,=\, \left\|\,y,\, b_{\,2},\, \cdots,\, b_{\,n}\,\right\|}\,\left|\,T\,(\,x,\, y,\, b_{\,2},\, \cdots,\, b_{\,n}\,)\,\right|$,
\item[$(iii)$]$\|\,T\,\| \,=\, \sup\limits_{\left\|\,x,\, b_{\,2},\, \cdots,\, b_{\,n}\,\right\| \,\neq\, 0,\; \left\|\,y,\, b_{\,2},\, \cdots,\, b_{\,n}\,\right\| \,\neq\, 0}\,\dfrac{\left|\,T\,(\,x,\, y,\, b_{\,2},\, \cdots,\, b_{\,n}\,)\,\right|}{\left\|\,x,\, b_{\,2},\, \cdots,\, b_{\,n}\,\right\|\,\left\|\,y,\, b_{\,2},\, \cdots,\, b_{\,n}\,\right\|}$,
\item[$(iv)$]$\left|\,T^{\,\prime}\,(\,x,\, b_{\,2},\, \cdots,\, b_{\,n}\,)\,\right| \,\leq\, \|\,T^{\,\prime}\,\|\,\left\|\,x,\, b_{\,2},\, \cdots,\, b_{\,n}\,\right\|^{\,2},\; \;\text{for all}\; x \,\in\, H$.
\item[$(v)$]$\|\,T^{\,\prime}\,\| \,=\, \sup\limits_{\left\|\,x,\, b_{\,2},\, \cdots,\, b_{\,n}\,\right\| \,=\, 1}\,\left|\,T^{\,\prime}\,(\,x,\, b_{\,2},\, \cdots,\, b_{\,n}\,)\,\right|$.
\end{description}
\end{remark}

\begin{theorem}\label{th2.2}
A b-sesquilinear functional \,$T$\, is bounded if and only if \,$T^{\,\prime}$\, is bounded.\,Moreover, \,$\|\,T^{\,\prime}\,\| \,\leq\, \|\,T\,\| \,\leq\, 2\,\|\,T^{\,\prime}\,\|$.
\end{theorem} 

\begin{proof}
First we suppose that \,$T$\, is bounded.\,Then for all \,$x \,\in\, H$, we have 
\begin{align*}
\left|\,T^{\,\prime}\,(\,x,\, b_{\,2},\, \cdots,\, b_{\,n}\,)\,\right| &\,=\, \left|\,T\,(\,x,\, x,\, b_{\,2},\, \cdots,\, b_{\,n}\,)\,\right| \\
&\,\leq\, \|\,T\,\|\,\left\|\,x,\, b_{\,2},\, \cdots,\, b_{\,n}\,\right\|\,\left\|\,x,\, b_{\,2},\, \cdots,\, b_{\,n}\,\right\|\\
& \,=\, \|\,T\,\|\,\left\|\,x,\, b_{\,2},\, \cdots,\, b_{\,n}\,\right\|^{\,2}. 
\end{align*}
So, \,$T^{\,\prime}$\, is bounded and \,$\|\,T^{\,\prime}\,\| \,\leq\, \|\,T\,\|$.\\

Conversely, suppose that \,$T^{\,\prime}$\, is bounded.\,By Theorem (\ref{lem2}) and using \,$(v)$\, of remark (\ref{note2}), for all \,$x,\,y \,\in\, H$, we obtain 
\begin{align*}
&\left|\,T\,(\,x,\, y,\, b_{\,2},\, \cdots,\, b_{\,n}\,)\,\right|\\
& \,\leq\, \dfrac{1}{4}\,\|\,T^{\,\prime}\,\|\,\left(\,\left\|\,x \,+\, y,\, b_{\,2},\, \cdots,\, b_{\,n}\,\right\|^{\,2} \,+\, \left\|\,x \,-\, y,\, b_{\,2},\, \cdots,\, b_{\,n}\,\right\|^{\,2}\,\right)\\
& \,+\, \dfrac{1}{4}\,\|\,T^{\,\prime}\,\|\,\left(\,\left\|\,x \,+\, i\,y,\, b_{\,2},\, \cdots,\, b_{\,n}\,\right\|^{\,2} \,+\, \left\|\,x \,-\, i\,y,\, b_{\,2},\, \cdots,\, b_{\,n}\,\right\|^{\,2}\,\right)\\
&=\, \dfrac{1}{4}\,\|\,T^{\,\prime}\,\|\,\left(\,4\,\left\|\,x,\, b_{\,2},\, \cdots,\, b_{\,n}\,\right\|^{\,2} \,+\, 4\,\left\|\,y,\, b_{\,2},\, \cdots,\, b_{\,n}\,\right\|^{\,2}\,\right)\; \;[\;\text{by Parallelogram law}\;]\\
&\,=\, \|\,T^{\,\prime}\,\|\,\left(\,\left\|\,x,\, b_{\,2},\, \cdots,\, b_{\,n}\,\right\|^{\,2} \,+\, \left\|\,y,\, b_{\,2},\, \cdots,\, b_{\,n}\,\right\|^{\,2}\,\right)\\
&\Rightarrow\, \sup\limits_{\left\|\,x,\, b_{\,2},\, \cdots,\, b_{\,n}\,\right\| \,=\, 1 \,=\, \left\|\,y,\, b_{\,2},\, \cdots,\, b_{\,n}\,\right\|}\,\left|\,T\,(\,x,\, y,\, b_{\,2},\, \cdots,\, b_{\,n}\,)\,\right| \,\leq\, 2\;\|\,T^{\,\prime}\,\|.    
\end{align*}
Therefore, \,$T$\, is bounded and \,$\|\,T\,\| \,\leq\, 2\,\|\,T^{\,\prime}\,\|$.\,This completes the proof.  
\end{proof} 

\begin{theorem}
If \,$T$\, is a bounded and symmetric \,$b$-sesquilinear functional then \,$\|\,T\,\| \,=\, \|\,T^{\,\prime}\,\|$.
\end{theorem} 

\begin{proof}
By Theorem (\ref{th2.2}), \,$T^{\,\prime}$\, is bounded and \,$\|\,T^{\,\prime}\,\| \,\leq\, \|\,T\,\|$.\,So, we need to prove that \,$\|\,T\,\| \,\leq\, \|\,T^{\,\prime}\,\|$.\,By Theorem (\ref{th2.1}), we note that \,$T^{\,\prime}$\, is real-valued.\,So from Theorem (\ref{lem2}), we obtain 
\begin{align*}
&\left|\,\text{Re}\,T\,(\,x,\, y,\, b_{\,2},\, \cdots,\, b_{\,n}\,)\,\right|\\
&\,\leq\, \dfrac{1}{4}\,\left(\,\left|\,T^{\,\prime}\left(\,x \,+\, y,\, b_{\,2},\, \cdots,\, b_{\,n}\,\right)\,\right| \,+\, \left|\,T^{\,\prime}\left(\,x \,-\, y,\, b_{\,2},\, \cdots,\, b_{\,n}\,\right)\,\right|\,\right).
\end{align*}
Using the boundedness of \,$T^{\,\prime}$\, and the Parallelogram law, we obtain
\begin{align*}
&\left|\,\text{Re}\,T\,(\,x,\, y,\, b_{\,2},\, \cdots,\, b_{\,n}\,)\,\right|\\
&\,\leq\, \dfrac{1}{4}\,\|\,T^{\,\prime}\,\|\,\left(\,\left\|\,x \,+\, y,\, b_{\,2},\, \cdots,\, b_{\,n}\,\right\|^{\,2} \,+\, \left\|\,x \,-\, y,\, b_{\,2},\, \cdots,\, b_{\,n}\,\right\|^{\,2}\,\right)\\
&=\, \dfrac{1}{4}\,\|\,T^{\,\prime}\,\|\,\left(\,2\,\left\|\,x,\, b_{\,2},\, \cdots,\, b_{\,n}\,\right\|^{\,2} \,+\, 2\,\left\|\,y,\, b_{\,2},\, \cdots,\, b_{\,n}\,\right\|^{\,2}\,\right).
\end{align*}
So, if \,$\left\|\,x,\, b_{\,2},\, \cdots,\, b_{\,n}\,\right\| \,=\, 1,\; \left\|\,y,\, b_{\,2},\, \cdots,\, b_{\,n}\,\right\| \,=\, 1$, we have
\begin{equation}\label{eq2.1}
\left|\,\text{Re}\,T\,(\,x,\, y,\, b_{\,2},\, \cdots,\, b_{\,n}\,)\,\right| \,\leq\, \|\,T^{\,\prime}\,\|.
\end{equation}
Writing \,$T\,(\,x,\, y,\, b_{\,2},\, \cdots,\, b_{\,n}\,)$\, in the form \,$r\,e^{\,i\,v}$\, and letting \,$\alpha \,=\, e^{\,-\, i\,v}$, we obtain
\begin{equation}\label{eq2.2}
\alpha\,T\,(\,x,\, y,\, b_{\,2},\, \cdots,\, b_{\,n}\,) \,=\, r \,=\, \left|\,T\,(\,x,\, y,\, b_{\,2},\, \cdots,\, b_{\,n}\,)\,\right|.
\end{equation}
Form (\ref{eq2.1}) and (\ref{eq2.2}),
\begin{align*}
\|\,T^{\,\prime}\,\|& \,\geq\, \left|\,\text{Re}\,T\,(\,\alpha\,x,\, y,\, b_{\,2},\, \cdots,\, b_{\,n}\,)\,\right| \,=\, \left|\,\text{Re}\,\alpha\,T\,(\,x,\, y,\, b_{\,2},\, \cdots,\, b_{\,n}\,)\,\right|\\
&=\, \left|\,T\,(\,x,\, y,\, b_{\,2},\, \cdots,\, b_{\,n}\,)\,\right|,
\end{align*}
whenever \,$\left\|\,x,\, b_{\,2},\, \cdots,\, b_{\,n}\,\right\| \,=\, 1,\; \left\|\,y,\, b_{\,2},\, \cdots,\, b_{\,n}\,\right\| \,=\, 1$.\,So,
\[\|\,T\,\| \,=\, \sup\limits_{\left\|\,x,\, b_{\,2},\, \cdots,\, b_{\,n}\,\right\| \,=\, 1 \,=\, \left\|\,y,\, b_{\,2},\, \cdots,\, b_{\,n}\,\right\|}\,\left|\,T\,(\,x,\, y,\, b_{\,2},\, \cdots,\, b_{\,n}\,)\,\right| \,\leq\, \|\,T^{\,\prime}\,\|.\hspace{2cm}\]
This proves the theorem.    
\end{proof} 
 
Now, we shall make use of the polarization identities to obtain a generalized form of the Schwarz inequality. 

\begin{theorem}(Generalized Schwarz inequality)
Let \,$T$\, be a positive bounded\;$b$-sesquilinear functional defined on \,$H \,\times\, H \,\times\,\left<\,b_{\,2}\,\right> \,\times\, \cdots \,\times\, \left<\,b_{\,n}\,\right>$\, and \,$T^{\,\prime}$\, be its associated quadratic form.\,Then
\[\left|\,T\,(\,x,\, y,\, b_{\,2},\, \cdots,\, b_{\,n}\,)\,\right|^{\,2} \,\leq\, T^{\,\prime}\,(\,x,\, b_{\,2},\, \cdots,\, b_{\,n}\,)\;T^{\,\prime}\,(\,y,\, b_{\,2},\, \cdots,\, b_{\,n}\,)\; \;\forall\; x,\,y \,\in\, H.\]
\end{theorem} 

\begin{proof}
We first note that if \,$T$\, is positive then we can write 
\[T^{\,\prime}\,(\,x,\, b_{\,2},\, \cdots,\, b_{\,n}\,) \,=\, T\,(\,x,\, x,\, b_{\,2},\, \cdots,\, b_{\,n}\,) \,\geq\, 0, \;\text{for all}\; x \,\in\, H\] i.\,e., \,$T^{\,\prime}$\, is real-valued and so by Theorem (\ref{th2.1}), \,$T$\, is symmetric, i.\,e., 
\[T\,(\,x,\, y,\, b_{\,2},\, \cdots,\, b_{\,n}\,) \,=\, \overline{\,T\,(\,y,\, x,\, b_{\,2},\, \cdots,\, b_{\,n}\,)},\;\text{for all}\; x,\,y \,\in\, H.\]If \,$T\,(\,x,\, y,\, b_{\,2},\, \cdots,\, b_{\,n}\,) \,=\, 0$\, then the inequality is clear.\,So we assume that \,$0 \,\neq\, T\,(\,x,\, y,\, b_{\,2},\, \cdots,\, b_{\,n}\,)$.\,For arbitrary scalars \,$\alpha,\, \beta$\, we obtain
\[0 \,\leq\, T^{\,\prime}\,(\,\alpha\,x \,+\, \beta\,y,\, b_{\,2},\, \cdots,\, b_{\,n}\,) \,=\, T\,(\,\alpha\,x \,+\, \beta\,y,\, \alpha\,x \,+\, \beta\,y,\, b_{\,2},\, \cdots,\, b_{\,n}\,)\]
\[=\, \alpha\,\overline{\alpha}\,T^{\,\prime}\,(\,x,\, b_{\,2},\, \cdots,\, b_{\,n}\,) \,+\, \alpha\,\overline{\beta}\,T\,(\,x,\, y,\, b_{\,2},\, \cdots,\, b_{\,n}\,)\hspace{2.7cm}\]
\[ \,+\, \overline{\alpha}\,\beta\,\,T\,(\,y,\, x,\, b_{\,2},\, \cdots,\, b_{\,n}\,) \,+\, \beta\,\overline{\beta}\,T^{\,\prime}\,(\,y,\, b_{\,2},\, \cdots,\, b_{\,n}\,)\hspace{2.7cm}\]
\[=\, \alpha\,\overline{\alpha}\,T^{\,\prime}\,(\,x,\, b_{\,2},\, \cdots,\, b_{\,n}\,) \,+\, \alpha\,\overline{\beta}\,T\,(\,x,\, y,\, b_{\,2},\, \cdots,\, b_{\,n}\,)\hspace{2.7cm}\]
\begin{equation}\label{eq2.3}
 \,+\, \overline{\alpha}\,\beta\,\,\overline{\,T\,(\,x,\, y,\, b_{\,2},\, \cdots,\, b_{\,n}\,)} \,+\, \beta\,\overline{\beta}\,T^{\,\prime}\,(\,y,\, b_{\,2},\, \cdots,\, b_{\,n}\,).\hspace{2.7cm}    
\end{equation}
Let \,$\alpha \,=\, t$\, be real and \,$\beta \,=\, \dfrac{T\,(\,x,\, y,\, b_{\,2},\, \cdots,\, b_{\,n}\,)}{\left|\,T\,(\,x,\, y,\, b_{\,2},\, \cdots,\, b_{\,n}\,)\,\right|}$.\,Then clearly, 
\[\overline{\beta}\,T\,(\,x,\, y,\, b_{\,2},\, \cdots,\, b_{\,n}\,) \,=\, \left|\,T\,(\,x,\, y,\, b_{\,2},\, \cdots,\, b_{\,n}\,)\,\right|\; \;\text{and}\; \;\beta\,\overline{\beta} \,=\, 1.\]
So, from (\ref{eq2.3}), for all real \,$t$, we have
\[0 \,\leq\, t^{\,2}\,T^{\,\prime}\,(\,x,\, b_{\,2},\, \cdots,\, b_{\,n}\,) \,+\, 2\,t\,\left|\,T\,(\,x,\, y,\, b_{\,2},\, \cdots,\, b_{\,n}\,)\,\right| \,+\, T^{\,\prime}\,(\,y,\, b_{\,2},\, \cdots,\, b_{\,n}\,).\]
So, the discriminant 
\[4\,\left|\,T\,(\,x,\, y,\, b_{\,2},\, \cdots,\, b_{\,n}\,)\,\right|^{\,2} \,-\, 4\,T^{\,\prime}\,(\,x,\, b_{\,2},\, \cdots,\, b_{\,n}\,)\;T^{\,\prime}\,(\,y,\, b_{\,2},\, \cdots,\, b_{\,n}\,)\]
cannot be positive.\,This proves the theorem.    
\end{proof} 

\begin{definition}
A linear operator \,$S \,:\, H \,\to\, H$\, is said to be \,$b$-bounded if there exists \,$M \,>\, 0$\, such that
\[\left\|\,S\,x,\, b_{\,2},\, \cdots,\, b_{\,n}\,\right\| \,\leq\, M\, \left\|\,x,\, b_{\,2},\, \cdots,\, b_{\,n}\,\right\|,\; \;\text{for all}\; \,x \,\in\, H.\]
The norm of \,$S$\, is defined as
\[\|\,S\,\| \,=\, \inf\,\left\{\,M \,>\, 0 \,:\, \left\|\,S\,x,\, b_{\,2},\, \cdots,\, b_{\,n}\,\right\| \,\leq\, M\, \left\|\,x,\, b_{\,2},\, \cdots,\, b_{\,n}\,\right\|\; \;\forall\; x \,\in\, H\,\right\}.\]
\end{definition}

By applying the Theorem (\ref{th2}), we finally give a general representation of \,$b$-sesquilinear functional in \,$n$-Hilbert space. 

\begin{theorem}
Let \,$T$\, be a bounded \,$b$-sesquilinear functional defined on \,$H \,\times\, H \,\times\,\left<\,b_{\,2}\,\right> \,\times\, \cdots \,\times\, \left<\,b_{\,n}\,\right>$\, and for each \,$x \,\in\, H$, the set \,$\left\{\,x,\, b_{\,2},\, \cdots,\, b_{\,n}\,\right\}$\, be linearly independent.\,Then \,$T$\, has a representation 
\begin{equation}\label{eq2.4}
T\,(\,x,\, y,\, b_{\,2},\, \cdots,\, b_{\,n}\,) \,=\, \left<\,S\,x,\, y \,|\, b_{\,2},\, \cdots,\, b_{\,n}\,\right>,\; \;\text{for all}\; x,\,y \,\in\, H,
\end{equation}
where \,$S \,:\, H \,\to\, H$\, is a \,$b$-bounded linear operator which is uniquely determined by \,$T$\, and \,$\|\,S\,\| \,=\, \|\,T\,\|$.
\end{theorem}

\begin{proof}
The functional \,$\overline{T\,(\,x,\, y,\, b_{\,2},\, \cdots,\, b_{\,n}\,)}$\, is linear in \,$y$\, because of the bar.\,Keeping \,$x$\, fixed, we apply Theorem (\ref{th2}), and obtain the representation
\[\overline{T\,(\,x,\, y,\, b_{\,2},\, \cdots,\, b_{\,n}\,)} \,=\, \left<\,y,\, z \,|\, b_{\,2},\, \cdots,\, b_{\,n}\,\right>,\]
where \,$y$\, is variable with \,$\left\{\,z,\, b_{\,2},\, \cdots,\, b_{\,n}\,\right\}$\, is linearly independent.\,So,
\begin{equation}\label{eq2.5}
T\,(\,x,\, y,\, b_{\,2},\, \cdots,\, b_{\,n}\,) \,=\, \left<\,z,\, y \,|\, b_{\,2},\, \cdots,\, b_{\,n}\,\right>. 
\end{equation}
Here \,$z \,\in\, H$\, is unique but depends on \,$x$.\,So we can write \,$z \,=\, S\,x$\, for some operator \,$S \,:\, H \,\to\, H$.\,In (\ref{eq2.5}), replacing \,$z$\, by \,$S\,x$, we obtain
\[T\,(\,x,\, y,\, b_{\,2},\, \cdots,\, b_{\,n}\,) \,=\, \left<\,S\,x,\, y \,|\, b_{\,2},\, \cdots,\, b_{\,n}\,\right>.\]
We now show that \,$S$\, is linear.\,For \,$\alpha,\, \beta \,\in\, \mathbb{K}$, we have
\begin{align*}
&\left<\,S\,\left(\,\alpha\,x_{\,1} \,+\, \beta\,x_{\,2}\,\right),\, y \,|\, b_{\,2},\, \cdots,\, b_{\,n}\,\right> \\
&\,=\, T\,(\,\alpha\,x_{\,1} \,+\, \beta\,x_{\,2},\, y,\, b_{\,2},\, \cdots,\, b_{\,n}\,)\\
&=\, \alpha\,T\,(\,x_{\,1},\, y,\, b_{\,2},\, \cdots,\, b_{\,n}\,) \,+\, \beta\,T\,(\,x_{\,2},\, y,\, b_{\,2},\, \cdots,\, b_{\,n}\,)\\
&=\, \alpha\,\left<\,S\,x_{\,1},\, y \,|\, b_{\,2},\, \cdots,\, b_{\,n}\,\right> \,+\, \beta\,\left<\,S\,x_{\,2},\, y \,|\, b_{\,2},\, \cdots,\, b_{\,n}\,\right>\\
& \,=\, \left<\,\alpha\,S\,x_{\,1} \,+\, \beta\,S\,x_{\,2},\, y \,|\, b_{\,2},\, \cdots,\, b_{\,n}\,\right>\; \;\forall\; y \,\in\, H.\\
&\Rightarrow\, S\,\left(\,\alpha\,x_{\,1} \,+\, \beta\,x_{\,2}\,\right) \,=\, \alpha\,S\,x_{\,1} \,+\, \beta\,S\,x_{\,2}.
\end{align*}
We now verify that \,$S$\, is \,$b$-bounded.\,If \,$S \,=\, 0$, there is nothing to prove.\,So, we assume that \,$S \,\neq\, 0$.\,From note (\ref{note2})\,$(iii)$ and equation (\ref{eq2.4}), we obtain
\begin{align*}
\|\,T\,\| &\,=\, \sup\limits_{\left\|\,x,\, b_{\,2},\, \cdots,\, b_{\,n}\,\right\| \,\neq\, 0,\; \left\|\,y,\, b_{\,2},\, \cdots,\, b_{\,n}\,\right\| \,\neq\, 0}\,\dfrac{\left|\,\left<\,S\,x,\, y \,|\, b_{\,2},\, \cdots,\, b_{\,n}\,\right>\,\right|}{\left\|\,x,\, b_{\,2},\, \cdots,\, b_{\,n}\,\right\|\,\left\|\,y,\, b_{\,2},\, \cdots,\, b_{\,n}\,\right\|}\\
&=\, \sup\limits_{\left\|\,x,\, b_{\,2},\, \cdots,\, b_{\,n}\,\right\| \,\neq\, 0,\; \left\|\,y,\, b_{\,2},\, \cdots,\, b_{\,n}\,\right\| \,\neq\, 0}\,\dfrac{\left|\,\left<\,S\,x,\, S\,x \,|\, b_{\,2},\, \cdots,\, b_{\,n}\,\right>\,\right|}{\left\|\,x,\, b_{\,2},\, \cdots,\, b_{\,n}\,\right\|\,\left\|\,S\,x,\, b_{\,2},\, \cdots,\, b_{\,n}\,\right\|}\\
&=\, \sup\limits_{\left\|\,x,\, b_{\,2},\, \cdots,\, b_{\,n}\,\right\| \,\neq\, 0}\,\dfrac{\left\|\,S\,x,\, b_{\,2},\, \cdots,\, b_{\,n}\,\right\|}{\left\|\,x,\, b_{\,2},\, \cdots,\, b_{\,n}\,\right\|}.
\end{align*} 
So,
\[\left\|\,S\,x,\, b_{\,2},\, \cdots,\, b_{\,n}\,\right\| \,\leq\, \,\|\,T\,\|\,\left\|\,x,\, b_{\,2},\, \cdots,\, b_{\,n}\,\right\|,\; \;\text{for all}\; x \,\in\, H.\]
Thus, it follows that \,$S$\, is \,$b$-bounded and moreover \,$\|\,S\,\| \,\leq\, \|\,T\,\|$.\,We have also by Schwarz inequality
\begin{align*}
\|\,T\,\| &\,=\, \sup\limits_{\left\|\,x,\, b_{\,2},\, \cdots,\, b_{\,n}\,\right\| \,\neq\, 0,\; \left\|\,y,\, b_{\,2},\, \cdots,\, b_{\,n}\,\right\| \,\neq\, 0}\,\dfrac{\left|\,\left<\,S\,x,\, y \,|\, b_{\,2},\, \cdots,\, b_{\,n}\,\right>\,\right|}{\left\|\,x,\, b_{\,2},\, \cdots,\, b_{\,n}\,\right\|\,\left\|\,y,\, b_{\,2},\, \cdots,\, b_{\,n}\,\right\|}\\
&\,\leq\, \sup\limits_{\left\|\,x,\, b_{\,2},\, \cdots,\, b_{\,n}\,\right\| \,\neq\, 0,\; \left\|\,y,\, b_{\,2},\, \cdots,\, b_{\,n}\,\right\| \,\neq\, 0}\,\dfrac{\left\|\,S\,x,\, b_{\,2},\, \cdots,\, b_{\,n}\,\right\|\,\left\|\,y,\, b_{\,2},\, \cdots,\, b_{\,n}\,\right\|}{\left\|\,x,\, b_{\,2},\, \cdots,\, b_{\,n}\,\right\|\,\left\|\,y,\, b_{\,2},\, \cdots,\, b_{\,n}\,\right\|} \,=\, \|\,S\,\|
\end{align*}
and so \,$\|\,S\,\| \,=\, \|\,T\,\|$.\,Now we show that \,$S$\, is unique.\,If possible suppose that there exists a linear operator \,$U \,:\, H \,\to\, H$\, such that 
\begin{align*}
&T\,(\,x,\, y,\, b_{\,2},\, \cdots,\, b_{\,n}\,) \,=\, \left<\,S\,x,\, y \,|\, b_{\,2},\, \cdots,\, b_{\,n}\,\right>\\
& \,=\, \left<\,U\,x,\, y \,|\, b_{\,2},\, \cdots,\, b_{\,n}\,\right>,\; \text{for all}\; \,x,\,y \,\in\, H.\\ 
&\Rightarrow\, \left<\,S\,x \,-\, U\,x,\, y \,|\, b_{\,2},\, \cdots,\, b_{\,n}\,\right> \,=\, 0,\; \;\text{for all}\; \,x,\,y \,\in\, H.
\end{align*}
In particular for \,$y \,=\, S\,x \,-\, U\,x$, 
\[\left<\,S\,x \,-\, U\,x,\, S\,x \,-\, U\,x  \,|\, b_{\,2},\, \cdots,\, b_{\,n}\,\right> \,=\, 0,\; \;\text{for all}\; x \,\in\, H.\]
\[\Rightarrow\, \left\|\,S\,x \,-\, U\,x,\, b_{\,2},\, \cdots,\, b_{\,n}\,\right\|^{\,2} \,=\, 0,\; \;\text{for all}\; x \,\in\, H.\hspace{1cm}\]
Since for each \,$x \,\in\, H$, the set \,$\left\{\,x,\, b_{\,2},\, \cdots,\, b_{\,n}\,\right\}$\, is linearly independent, \,$S\,x \,=\, U\,x$\, for all \,$x$\, and so \,$S \,=\, U$.\,This proves the theorem.
\end{proof}

\end{document}